\documentclass{elsarticle}
\usepackage{graphicx}
\usepackage{amsmath}
\usepackage{amssymb}
\usepackage{caption}
\usepackage{stmaryrd}
\usepackage{multirow}

\newcommand{\con}{\mathrm{fer}}
\newcommand{\ub}{\mathbf{u}}
\newcommand{\adj}{\xi}
\newcommand{\iso}{\mathrm{iso}}

\newcommand{\nud}{\boldsymbol{\nu}_{\mathrm{d}}}
\newcommand{\nub}{\boldsymbol{\nu}}
\newcommand{\nul}{\boldsymbol{\nu}_{\mathrm{L}}}
\newcommand{\err}{\varepsilon}
\newcommand{\pol}{\mathrm{pol}}
\newcommand{\loc}{\mathrm{loc}}
\DeclareMathOperator{\rot}{\nabla \times}
\let\div\relax
\DeclareMathOperator{\div}{\nabla \cdot}
\usepackage{color}

\newdefinition{rmk}{Remark}

\newdefinition{as}{Assumption}

\usepackage{lineno,hyperref}
\modulolinenumbers[5]

\begin{document}

\begin{frontmatter}

\title{A Defect Corrected Finite Element Approach for the Accurate Evaluation of Magnetic Fields on Unstructured Grids}

\author[mymainaddress,mysecondaryaddress]{Ulrich R\"omer\corref{mycorrespondingauthor}}
\cortext[mycorrespondingauthor]{Corresponding author}
\ead{roemer@temf.tu.darmstadt.de}

\author[mymainaddress,mysecondaryaddress]{Sebastian Sch\"ops}

\author[mymainaddress,mysecondaryaddress]{Herbert De Gersem}

\address[mymainaddress]{Institut f\"ur Theorie Elektromagnetischer Felder, Technische Universit\"at Darmstadt, Schlossgartenstr. 8, D-64289 Darmstadt, Germany}
\address[mysecondaryaddress]{Graduate School of Computational Engineering, Technische Universit\"at Darmstadt, Dolivostr 15, D-64293 Darmstadt, Germany}

\begin{abstract}
In electromagnetic simulations of magnets and machines one is often interested in a highly accurate and local evaluation of the magnetic field uniformity. Based on local post-processing of the solution, a defect correction scheme is proposed as an easy to realize alternative to higher order finite element or hybrid approaches. Radial basis functions (RBF)s are key for the generality of the method, which in particular can handle unstructured grids. Also, contrary to conventional finite element basis functions, higher derivatives of the solution can be evaluated, as required, e.g., for deflection magnets. Defect correction is applied to obtain a solution with improved accuracy and adjoint techniques are used to estimate the remaining error for a specific quantity of interest. Significantly improved (local) convergence orders are obtained. The scheme is also applied to the simulation of a Stern-Gerlach magnet currently in operation.
\end{abstract}


\end{frontmatter}

\section{Introduction}
The modern design of complex electromagnetic devices is based on efficient and accurate higher-order computational schemes. Despite significant improvements, specific challenges concerning higher-order modeling persist for a large variety of physical models. Frequently in computational magnetics one is interested in local rather than global (energy-related) quantities of interest. Important application examples are electrical machines and magnets used for particle deflection or focusing in accelerating structures \cite{Russenschuck2011,salon1995finite}. The design goal is a locally uniform magnetic field in an air gap, quantified, e.g., by Fourier harmonics. These Fourier harmonics are a key input for particle tracking codes and are used for electromotive force or cogging torque computations of electrical machines. Assessing the field uniformity can be a difficult task due to the following: the air gap constitutes only a small fraction of the computational domain, mainly consisting of large iron or steel parts with a complicated geometry. Also, Fourier harmonics need to be evaluated with high accuracy, as in the case of accelerator magnets, higher order coefficients of small magnitude can cause beam instabilities. Moreover, in deflection magnets the quantity of interest is a uni-directional derivative of the magnetic flux density, which is typically not well-defined in a finite element (FE) approach involving vector potentials.    

Dedicated schemes have been presented in the literature, exploiting the smoothness of the solution in the air gap: a boundary element - finite element coupling \cite{Kurz1998} and a hybrid finite element - spectral element scheme \cite{DeGersem2002}, among others. However, both approaches require significant code modifications. As mentioned above, the simulation of deflection magnet applications is particularly challenging. Even in combination with a higher order finite element approach \cite{demkowicz2006computing,dular2009posteriori}, a dedicated local post-processing is required. Such a solution reconstruction has been presented in \cite{DeGersem2014} based on an analytical solution. A higher differentiability of the solution across the element boundaries could also be guaranteed by using isogeometric finite elements as shown in \cite{pelsoptimization}. In isogeometric methods spline basis functions of arbitrary regularity can be used. However, local refinement strategies are still difficult to realize and further research is needed make applications to complex geometries possible in an automated way. 

In this work another approach is investigated, referred to as defect correction in the literature \cite{Stetter1978,bohmer2012defect,Giles2002}. It is based on a solution reconstruction, however, additionally a part of the numerical error is estimated and removed from the solution to obtain a faster convergence. Adjoint correction \cite{pierce2000adjoint,Giles2002} is applied to estimate the remaining error for a quantity of interest. Following the approach outlined in \cite{Pierce2004}, no additional higher-order discrete operator is required. Defect correction schemes are often applied on structured grids based on a tensor-product spline reconstruction, see \cite{Giles2002}. The case of an unstructured grid did not receive much attention so far. A velocity reconstruction of fluid flows on unstructured grids was discussed in \cite{basumatary2014defect}. In \cite{Giles2002} biharmonic smoothing was presented and analyzed. As solvers for the biharmonic equation are typically unavailable in a computational magnetics context, a generic post-processing by radial basis functions (RBFs) is presented as an alternative in this paper. RBFs have been already successfully used to post-process solutions for scalar and vector fields \cite{bonaventura2011kernel}. A key contribution here is the discussion of a local defect correction approach and the associated numerical errors. Numerical results are given to illustrate the convergence orders and the accuracy of the defect corrected quantities of interest using RBFs. Although motivated from an electromagnetic perspective, post-processing by RBFs and the scheme in general is not limited to magnetic field problems.

The paper is structured as follows: in Section \ref{sec:magnetics} a simplified magnetostatic model problem is formulated together with a FE scheme. Defect correction principles are presented in Section \ref{sec:dec}, together with RBF post-processing and convergence results. In Section \ref{sec:interface} the model is extended to interface problems with more complicated geometries, to cover practical applications. Defect correction is also adapted to this more general setting. Quantities of interest, appearing in practical applications are discussed in Section \ref{sec:qoi}. The findings are illustrated by an academic example and the simulation of a Stern-Gerlach magnet, currently in operation at KU Leuven, in Section \ref{sec:numerics}.

\section{Newton Method for Magnetostatics}
\label{sec:magnetics}
In a classical setting, macroscopic magnetic fields are governed by Maxwell's equations. Solving the full Maxwell system can be challenging and is often unnecessary for devices operated at low-frequencies, or even in a stationary regime. In particular, in the static limit we obtain the magnetostatic problem 
\begin{subequations}
\begin{align}
\div \vec{B} &= 0,  \quad \ \mathrm{in} \ \Omega,\\
\rot \vec{H} &= \vec{J},  \quad \ \mathrm{in} \ \Omega,\\
\vec{B} \cdot \vec{n} &= 0, \quad \ \mathrm{on} \ \partial \Omega.
\end{align}
\label{eq:ms}
\end{subequations}
In \eqref{eq:ms}, $\Omega$ denotes the (bounded) computational domain and $\vec{B},\vec{H}$ refer to the magnetic flux density and magnetic field strength, respectively. Also, $\vec{n}$ and $\vec{J}$ denote the outer unit normal and the source current density, respectively. The material constitutive relation reads $\vec{H}=\nu(|\vec{B}|) \vec{B}$, where $\nu(|\vec{B}|)$ refers to the magnetic reluctivity, which may depend on the magnitude of the magnetic flux density due to ferromagnetic saturation \cite{cullity2011introduction}. 

Introducing the magnetic vector potential $\vec{B} = \rot \vec{A}$ in \eqref{eq:ms}, we are concerned with the second order system
\begin{subequations}
\begin{align}
\rot \left (\nu(|\rot \vec{A}|) \rot \vec{A}\right) &= \vec{J},  \quad \ \mathrm{in} \ \Omega,\\
\vec{A} \times \vec{n} &= 0, \quad \ \mathrm{on} \ \partial \Omega.
\end{align}
\label{eq:ms_so}
\end{subequations}
For many applications, or in an early design phase, a two-dimensional setup can be considered. This is the case, e.g., for a geometry invariant with respect to translations in one coordinate direction. Then, \eqref{eq:ms_so} reduces to 
\begin{subequations}
\begin{align}
-\nabla \cdot \left( \nu(|\nabla u|) \nabla u \right) &= J \quad \ \mathrm{in} \ \Omega,\\
u &= 0, \quad \ \mathrm{on} \ \partial \Omega,
\end{align}
\label{eq:ms_td}
\end{subequations}
where $u$ and $J$ refer to the remaining components of the magnetic vector potential and source current density, respectively. For a nonlinear problem such as \eqref{eq:ms_td}, defect correction can be realized as an additional (approximate) Newton iteration \cite{Pierce2004}. Hence, we linearize \eqref{eq:ms_td} as 
\begin{subequations}
\begin{align}
-\nabla \cdot \left (\nul(\nabla u^{(k-1)}) \nabla u^{(k)} \right )&= j(\nabla u^{(k-1)}), \quad \ \mathrm{in} \ \Omega,\\
u &= 0, \quad \ \mathrm{on} \ \partial \Omega,
\end{align} 
\label{eq:ms_lin}
\end{subequations}
at step $k$, with the tensor 
\begin{equation}
\nul(\vec{r}) = \nu(|\vec{r}|) \mathbb{I} + 
\left \{
\begin{aligned}
&\frac{\nu^{'}(|\vec{r}|)}{|\vec{r}|} \vec{r} \otimes \vec{r},  \quad |\vec{r}| \neq 0, \\
&0, \quad |\vec{r}| = 0, 
\end{aligned}
\right .
\end{equation}
where $\nu^{'}(x) := \mathrm{d} \nu(x)/ \mathrm{d}x$, $\mathbb{I}$ refers to the $2 \times 2$ identity matrix and $\vec{r}$ is an arbitrary vector. Setting $j(\vec{r}) = j_\mathrm{N}(\vec{r}) = J + \nabla \cdot \left(\nu(|\vec{r}|)\vec{r} -\nul(\vec{r})\vec{r}\right)$ we obtain the Newton-Raphson method. However, we allow for more general sources $j$. In the remaining part of the paper, the index $k$ and the subscript $\mathrm{L}$ are omitted for simplicity.
\begin{rmk}
Note that in this two-dimensional setting, $\nul$ is related to the differential reluctivity tensor $\nud(\vec{r}) = D_{\vec{r}} \nu(|\vec{r}|)$ as 
\begin{equation}
\nul = 
\begin{pmatrix}
\nu_{d,22} & -\nu_{d,21} \\
-\nu_{d,12} & \nu_{d,11}
\end{pmatrix} 
.
\end{equation}
\end{rmk}
To \eqref{eq:ms_lin} we associate the weak formulation, find $u \in V = H_0^1(\Omega)$ such that
\begin{equation}
\int_{\Omega} \nub \nabla u \cdot \nabla v \ \mathrm{d} x = \int_{\Omega} j v \ \mathrm{d} x, 
\label{eq:msweak}
\end{equation}
for all $v \in V$. Equation (\ref{eq:msweak}) is discretized by the FE method on a triangular mesh. Let $(\vec{x}_i)_{i=1}^N$ be the nodes of the mesh, $\mathcal{T}_h$ the set of elements and $\mathcal{P}^k(K)$ denote the space of polynomials with degree $p \leq k$. We seek $u_h$ in 
\begin{equation}
V_h = \{ v_h \in \mathcal{C}(\Omega) \ | \ v_h |_K \in \mathcal{P}^1(K) \ \forall K \in \mathcal{T}_h, v_h = 0 \ \mathrm{on} \ \partial \Omega\}.
\end{equation}
We also have $V_h = \mathrm{span}\{ \phi_i, i=1,\dots,N_{\mathrm{int}} \}$, where $\phi_i$ denote piecewise linear and continuous shape functions associated to interior nodes. The degrees of freedom $\ub$ are subject to the linear system of equations 
\begin{equation}
\mathbb{K} \ub = \mathbf{j}, \quad \mathbb{K}_{ij} = \int_{\Omega} \nub \nabla \phi_j \cdot \nabla \phi_i \ \mathrm{d} x, \quad j_i = \int_{\Omega} j \phi_i \ \mathrm{d} x, \quad i,j=1,\dots,N_\mathrm{int}.
\label{eq:}
\end{equation}
Note that there is no conceptual difference with respect to discretization in the nonlinear setting. 

In applications one is interested not in the solution $u$ itself, but in a quantity of interest $\mathcal{F}$. In computational magnetics applications, $\mathcal{F}$ typically refers to local measures of the field uniformity, such as Fourier harmonics of the solution, or magnetic forces. A detailed discussion of quantities of interest is given in Section \ref{sec:qoi}. 

\section{Defect Correction}
\label{sec:dec}
The fundamental idea of defect correction is to interpolate the low-order numerical solution between the nodes of the mesh to obtain a higher order reconstruction. This reconstruction can be used to both estimate and reduce the error by solving the finite element equation an additional time. Different versions of defect correction schemes have been thoroughly addressed in the literature for the finite difference and finite element method and we refer to \cite{Stetter1978,Pierce2000,Giles2002,Pierce2004} and the references therein. One can either improve the accuracy of the solution directly, referred to as primal approach here, or use the adjoint solution to obtain an improved quantity of interest. Both approaches can be combined to achieve an ever better rate of convergence. However, in this work the adjoint approach is used to estimate the error after primal correction, solely. 

To simplify the exposition of defect correction principles we assume for the time-being:
\begin{as}
\label{as:smooth}
The solution $u$ is smooth. 
\end{as}
Of course, the smoothness of the solution does not hold true for any realistic setting and the assumption will be relaxed later on.

\subsection{Primal Approach}
Let $\pi_h u_h$ be a reconstruction of the solution, where $\pi_h$ is an operator with approximation accuracy
\begin{equation}
\|v - \pi_h v\|_{L^2(\Omega)} = \mathcal{O}(h^r),
\label{eq:order_rec}
\end{equation}
for $2 \leq r \leq 4$. Following \cite{Giles2002}, in the primal approach we solve for the correction $e_h \in V_h$ subject to
\begin{equation}
\int_{\Omega} \nub \nabla e_h \cdot \nabla v_h \ \mathrm{d} x = \int_{\Omega} j v_h \ \mathrm{d} x - \int_{\Omega} \nub \nabla \pi_h u_h \cdot \nabla v_h \ \mathrm{d} x, \quad \forall v_h \in V_h.
\label{eq:dec_prim}
\end{equation}
Then an improved solution is obtained as $\tilde{u}_h := \pi_h (u_h + e_h)$. We recall from \cite{Pierce2004}, that the present defect correction approach only requires the assembly of a new right-hand-side, which is in contrast to early defect correction procedures, based on higher-order discrete operators. This assembly, however, demands for a higher-order numerical quadrature. 

The decay rate of the remaining error in $\tilde{u}_h$ in the $L^2$-norm is now given by the reconstruction accuracy $r$ \cite{Pierce2004}. Also, the convergence rate of the output functional $\mathcal{F}$ exhibits the same improvement. Defect correction can be applied repeatedly yielding a further improved solution, however, with a convergence rate still limited to $r$ \cite{barrett1988optimal}. 
As an example, for a bivariate cubic $\mathcal{C}^2$ spline reconstruction on a structured grid, the error after defect correction decays as $\mathcal{O}(h^4)$ \cite{Giles2002}. Here, we allow for an unstructured grid and tensor product reconstructions are not applicable if interpolation is carried out at the nodes. This however, is crucial for the accuracy of the reconstruction, see \cite{barrett1988optimal}. To this end we propose the use of Radial Basis Functions (RBF) in Section \ref{sec:rbf}.

\subsection{Adjoint Approach}
Adjoint techniques can be used to estimate the error in the improved quantity of interest $\mathcal{F}(\tilde{u}_h)$ \cite{pierce2000adjoint}. For simplicity, we consider a linear quantity of interest in this section, i.e., $\mathcal{F}(\tilde{u}_h) = (g,\tilde{u}_h)$, where $(\cdot,\cdot)$ refers to the $L^2$-inner product. The adjoint solution $\adj \in V$ is subject to
\begin{subequations}
\begin{align}
-\nabla \cdot \left (\nub^\top \nabla \adj \right )&= g, \quad \ \mathrm{in} \ \Omega,\\
\adj &= 0, \quad \ \mathrm{on} \ \partial \Omega.
\end{align} 
\label{eq:ms_adjoint}
\end{subequations}
Problem \eqref{eq:ms_adjoint} can be approximated, again using the finite element method on the same grid with the same polynomial approximation functions. Let $\adj_h$ denote the associated discrete adjoint variable. Using the properties of the adjoint operator we infer
\begin{align}
\mathcal{F}(u - \tilde{u}_h) &= (g,u - \tilde{u}_h) = (\adj,j + \nabla \cdot  (\nub \nabla  \tilde{u}_h)) \notag\\
& = (\pi_h \adj_h,j + \nabla \cdot  (\nub \nabla  \tilde{u}_h)) - (\pi_h \adj_h - \adj,j + \nabla \cdot (\nub \nabla \tilde{u}_h)),
\label{eq:adjoint_error}
\end{align}
cf. \cite{Pierce2004}. As the second term on the right-hand-side of the previous expression is of higher-order, it can be neglected and hence, 
\begin{equation}
(\pi_h \adj_h,j + \nabla \cdot  (\nub \nabla  \tilde{u}_h))
\label{eq:dum}
\end{equation}
provides an asymptotically exact error bound.
Note that this error could also be removed to obtain an even higher order of convergence $\mathcal{O}(h^{r} + h^{\min(2,r-2)})$, see \cite{Pierce2004}.

\subsection{Radial Basis Functions}
\label{sec:rbf}
RBF interpolation is a widely used technique for the interpolation of scattered data, see \cite{buhmann_radial_2000} for an overview. We determine  
\begin{equation}
\pi_h u_h(\vec{x}) = \sum_{i=1}^N \alpha_i \Phi(|\vec{x}-\vec{x}_i|) + p(\vec{x}),
\label{eq:rbf}
\end{equation}
where $|\cdot|$ refers to the Euclidean norm, such that $\pi_h u_h(\vec{x}_i) = u_h(\vec{x}_i)$ and $\sum_{i=1}^N \alpha_i q(\vec{x}_i) = 0$ for all $q \in \mathcal{P}^{m-1}$, the space of (global) polynomials of degree less than $m$. In \eqref{eq:rbf}, the polynomial $p$ is required to ensure existence and uniqueness of $\pi_h u_h$, depending on the type of the RBF used. In particular we consider the following instances of polyharmonic splines, see \cite{buhmann_radial_2000},
\begin{equation}
\Phi_k(|\vec{x}-\vec{x}_i|) = 
\left \{
\begin{aligned}
&|\vec{x}-\vec{x}_i|^{2} \log(|\vec{x}-\vec{x}_i|), \ k=1,\\
&|\vec{x}-\vec{x}_i|^{3}, \  k=2, \\
&|\vec{x}-\vec{x}_i|^{5}, \  k=3 \\
\end{aligned}
\right . .
\label{eq:thin_plate}
\end{equation}
For $k=1,2$ in \eqref{eq:thin_plate} we have $m = 2$, whereas for $k=3$, $m=3$ holds. In several cases, the restrictions to be imposed on the nodes in order to render the interpolation problem well-posed, can be easily verified. Indeed, for $m=2$ we require pairwise distinct nodes that do not form a subset of a straight line \cite[p.3]{buhmann_radial_2000}. This is clearly fulfilled for a finite element mesh. 
Introducing a basis $(p_1,\dots,p_{M_m})$ of $\mathcal{P}^{m-1}$, the interpolation problem reads
\begin{equation}			
\begin{bmatrix}
\mathbb{F}_k & \mathbb{P}^\top \\
\mathbb{P} & 0
\end{bmatrix}
\begin{bmatrix}
\mathbf{a} \\ \mathbf{b}
\end{bmatrix}
=
\begin{bmatrix}
\ub \\ 0
\end{bmatrix}
,
\label{eq:rbf_system}
\end{equation}
where $(\mathbb{F}_k)_{ij} = \Phi_k (|\vec{x}_i - \vec{x}_j|)$, $\mathbb{P}_{ji} = p_j(\vec{x}_i)$, and $a_i=\alpha_i$, cf. \cite{gumerov2007fast}. Solving \eqref{eq:rbf_system} demands for dedicated numerical schemes, as dense matrices are involved. This is due to the non-local support of radial basis functions. In \cite{gumerov2007fast} a preconditioned Krylov method with $\mathcal{O}(N \log{N})$ complexity has been proposed. It is based on an acceleration of matrix-vector products based on the fast multipole method.

It should be noted that the RBF reconstruction is non-conforming, as $\pi_h u_h$ does not have a vanishing trace and hence $\pi_h u_n \notin V$. A theoretical investigation of this effect remains open. In all numerical experiments reported in this paper, this boundary residual did not have any noteworthy influence. Also, the additional term in the adjoint error estimate \eqref{eq:adjoint_error} 
\begin{equation}
\int_{\partial \Omega} \vec{n} \cdot (\nub^\top \nabla \xi)  \tilde{u}_h \ \mathrm{d} x,
\label{eq:boundary_contribution}
\end{equation}
which can be approximated by replacing $\xi$ by its discrete counterpart, was found to be negligible.  

Neglecting the aforementioned errors at the boundary, the convergence order of the defect correction method is identical to the approximation order of the RBFs. Hence, we conclude the section by recalling these estimates from the literature \cite{johnson2001}. Let $\pi_k$ denote the RBF reconstruction based on $\Phi_k$, $k=1,2,3$. We have the following approximation orders
\begin{equation}
\| v - \pi_k v \|_{L^2(\Omega)} = 
\left \{  
\begin{aligned}
&\mathcal{O}(h^{5/2}), \ k=1, \\
&\mathcal{O}(h^{3}), \ k=2, \\
&\mathcal{O}(h^{4}), \ k=3, \\
\end{aligned}
\right .
\label{eq:rbf_order}
\end{equation}
provided that the boundary is Lipschitz continuous. These results can be improved if the support of $v$ is a compact subset of $\Omega$ \cite[Theorem 3]{buhmann_radial_2000}. Note that the approximation orders for $k=2,3$ in \eqref{eq:rbf_order} were derived in \cite{johnson2001} for the case $\mathbb{R}^d$, with $d$ odd, solely. However, they seem to hold for $d$ even too, as they have been partially used in this way in \cite{buhmann_radial_2000}.

\begin{rmk}
In the RBF context, approximation orders are typically expressed by the global data density \cite{wendland_meshless_1999}
\begin{equation}
\tilde{h} := \sup \limits_{\vec{x} \in \Omega} \min \limits_{i=1,\dots,N} |\vec{x} - \vec{x}_i|.
\label{eq:h_rbf}
\end{equation}
The results remain true, however, for the mesh size $h$, as $\tilde{h} \leq h$.
\end{rmk}

\begin{rmk}
Additional efforts are needed to address the three-dimensional case with vector fields in $H(\mathrm{curl},\Omega)$. Vector-valued radial basis functions have been proposed in \cite{fuselier2007refined}. In \cite{bonaventura2011kernel} a vectorial reconstruction was presented, for the case of $H(\mathrm{div},\Omega)$ vector fields, improving finite element approximations with Raviart-Thomas elements but without defect correction.
\end{rmk} 

\begin{rmk}
In this work, only time-independent problems are addressed. However, defect correction can also be applied to initial boundary value problems. For instance, in \cite{giles2004progress} Burgers' equation has been considered with a solution reconstruction in space and time. Both primal and adjoint techniques are applicable in this case. However, for time-dependent nonlinear problems, solving the adjoint problem can be very challenging. This is due to the fact that the adjoint problem is solved in reverse time-direction requiring the primal solution at every discrete point in time for linearization. 
\end{rmk}

\section{Interface Problem}
\label{sec:interface}
\begin{figure}[t!]
\begin{minipage}[t!]{0.5\textwidth}
\includegraphics[width=0.9\textwidth]{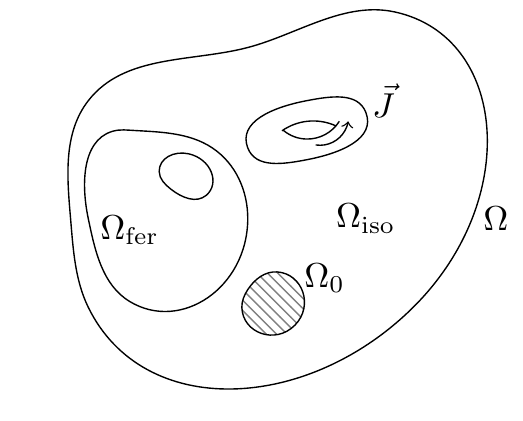}
\end{minipage}
\begin{minipage}[t!]{0.5\textwidth}
\includegraphics[width=0.9\textwidth]{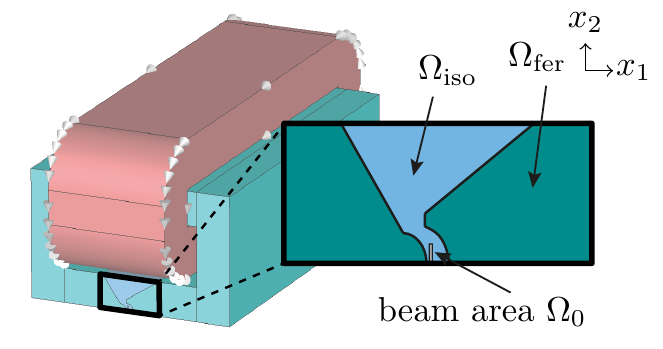}
\end{minipage}
\caption{Sketch of the model geometry. Computational domain with ferromagnetic and non-ferromagnetic domain $\Omega_\con$ and $\Omega_\iso$, respectively. Quantities of interest are evaluated in $\Omega_0$, solely.}
\label{fig:intro}
\end{figure}

We reformulate the geometrical setup in the following, to address applications with more complex and less regular geometries. A sketch of the setup we have in mind is given in Figure \ref{fig:intro}. In particular, the computational domain $\Omega$ is decomposed into a ferromagnetic and non-ferromagnetic domain $\Omega_\con$ and $\Omega_\iso = \Omega \setminus \overline{\Omega_\con}$, respectively. This is the typical setting of an interface problem. The magnetic reluctivity tensor is assumed to be constant $\nub(\vec{x},\cdot) = \nu_0$, for $\vec{x} \in \Omega_\iso$, but inhomogeneous in $\Omega_\con$. For simplicity we omit the explicit dependency of $\nub$ on $\vec{x}$.
Magnetic fields, and in particular field uniformity is evaluated in $\Omega_0 \subset \Omega_\iso$. 

\subsection{Local Defect Correction}
A serious difficulty arises when applying the defect correction approach to interface problems. To see this, we introduce a jump operator, e.g., at the material interface as
\begin{equation}
\llbracket \mathbf{n} \cdot (\nub \nabla u) \rrbracket = \mathbf{n}_{\con} \cdot (\nub \nabla u)_{\con} + \mathbf{n}_{\iso} \cdot (\nub \nabla u)_{\iso},
\end{equation}
where the subscripts $\con$ and $\iso$ denote restrictions to both sides of the interface $\Omega_\con$ and $\Omega_\iso$, respectively. Then, the solution $u$ fulfills the interface condition 
\begin{equation}
\llbracket  \mathbf{n} \cdot \left (\nub \nabla u \right )\rrbracket = 0, \ \mathrm{on} \ \partial \Omega_\con.
\label{eq:int}
\end{equation}  
Unfortunately, it is very difficult to define an RBF reconstruction $\pi_k u_h$ such that \eqref{eq:int} remains valid. Moreover, most interface configurations in practice contain geometrical singularities. Hence, Assumption \ref{as:smooth} does not hold, i.e., the solution is not smooth. We therefore localize the defect correction scheme to the area, where the quantities of interest are finally evaluated. Local defect correction schemes have been considered before, in particular we mention \cite{hackbusch1984local}. In contrast to the present approach, typically structured grids with a locally refined grid are used and several iterations between coarse and fine grid are carried out. 

As we are interested in the error in $\Omega_0$ solely, we can restrict ourselves to a local reconstruction. The discussion is based on the following assumption: 
\begin{as}
The domain $\Omega_0$ is simply connected with a polygonal Lipschitz boundary. The source current vanishes in $\Omega_0$.
\label{as:smooth_loc}
\end{as}
We observe that, in $\Omega_0$, the solution $u$ is smooth and also that the RBF approximation orders hold true, due to the previous assumption. Moreover, we can neglect boundary approximation errors, as $\partial \Omega_0$ is polygonal. Let $\pi_{h,0}$ denote the associated local RBF reconstruction operator and 
\begin{equation}
V_{h,0} = \{ v_{h,0} = v_h |_{\Omega_0}, v_h \in V_h \  |  \  v_{h,0}=0, \ \mathrm{on} \ \partial \Omega_0\}.
\end{equation}
Starting from \eqref{eq:dec_prim}, again dropping the iteration index, we solve for the correction $e_{h,0} \in V_{h,0}$ subject to
\begin{equation}
\int_{\Omega_0} \nu_0 \nabla e_{h,0} \cdot \nabla v_{h,0} \ \mathrm{d} x = - \int_{\Omega_0} \nu_0 \nabla \pi_{h,0} u_{h,0} \cdot \nabla v_{h,0} \ \mathrm{d} x, \quad \forall v_{h,0} \in V_{h,0}.
\label{eq:dec_prim_loc}
\end{equation}
Then an improved solution is obtained as $\tilde{u}_{h,0} := \pi_{h,0}(u_{h,0} + e_{h,0})$. Note that the local reconstruction $\pi_{h,0} u_{h,0}$ might be computed using a direct solver with a cost of $\mathcal{O}(N_0^3)$ as the number of FE nodes in $\Omega_0$ is small. 

By doing so, we correct the local approximation error, solely. Errors arising, e.g., at interface singularities will pollute into $\Omega_0$ unless local mesh refinement has been applied. To see this, we perform an error splitting on $\Omega_0$ as {$\err_0 := u - u_{h,0} = \err_0^{\loc} + \err_0^{\pol}$}. Then, the local error $\err_0^{\loc}$ is given as the FE error of the problem
\begin{align}
\nabla \cdot (\nu_0 \nabla u_0^{\loc} ) &= 0, \ \mathrm{in} \ \Omega_0, \\
u_0^{\loc} &= u, \ \mathrm{on} \ \partial \Omega_0,
\end{align}
whereas the pollution error $\err_0^{\pol}$ is the finite element error of the problem
\begin{align}
\nabla \cdot (\nu_0 \nabla u_0^{\pol} ) &= 0, \ \mathrm{in} \ \Omega_0, \\
u_0^{\pol} &= u - u_h, \ \mathrm{on} \ \partial \Omega_0,
\end{align}
cf. \cite[Remark 2.1.2]{babuvska1997pollution}. Let $V_0= H_0^{1}(\Omega_0)$, then the local error $\err_0^{\loc}$ is subject to 
\begin{equation}
\int_{\Omega_0} \nu_0 \nabla \err_0^{\loc} \cdot \nabla v_0 \ \mathrm{d} x = -\int_{\Omega_0} \nu_0 \nabla u_{h,0} \cdot \nabla v_0 \ \mathrm{d} x, \ \forall v_0 \in V_0.
\label{eq:loc_error}
\end{equation}
Note that $\err_0^{\loc}$ cannot simply be approximated by restricting $V_0$ to $V_{h,0}$ in \eqref{eq:loc_error}, as the right-hand-side would vanish. However, using the reconstructed solution, \eqref{eq:dec_prim_loc} represents a suitable approximation to \eqref{eq:loc_error}. Now the question arises, how the pollution error can be controlled locally, i.e., in $\Omega_0$. Dedicated techniques have been developed to this end, see \cite{oden1996local,babuvska1997pollution}. Here, we simply control the error in the global $H^1$-norm using an explicit residual error estimator. Following \cite{ainsworth1997posteriori},
\begin{equation}
\eta = \left(\sum_{K \in \mathcal{T}_h} h_K^2 \|j + \nabla \cdot (\nub \nabla u_h)\|_{L^2(K)}^2 + \sum_{\gamma \in \partial K} h_K \| \frac{1}{2} \llbracket (\nub \nabla u_h) \cdot \mathbf{n} \rrbracket \|_{L^2(\gamma)}^2 \right)^{1/2}, 
\label{eq:error_est}
\end{equation}
represents a reliable and efficient estimator for the FE error $\|u - u_h\|_{H^1_0(\Omega)}$. In \eqref{eq:error_est} local error contributions $\eta_K$ can also be identified, which are useful for an adaptive refinement process.

Provided that the pollution error is small enough, the defect correction scheme yields the same convergence orders as in the beginning of this section under the assumption of a smooth solution.

\section{Quantities of Interest}
\label{sec:qoi}
Typically, $\Omega_0$ is a circular or rectangular domain with a center identical to the center of the $x_1-x_2$ plane, see Figure \ref{fig:intro}. In the following, we consider two particular examples of quantities of interest $\mathcal{F}: H^m(\Omega_0) \rightarrow \mathbb{R}$, where we allow for $m > 1$. For magnet and machines applications we compute Fourier coefficients of the solution. Note that in the accelerator literature the notion of multipole coefficients is usually preferred \cite{Russenschuck2011}. For deflection magnets applications one is rather interested in (averaged) uni-directional derivatives of the magnetic flux density.

\subsection{Fourier Coefficients}
Fourier coefficients are typically extracted from $u$ in the beam pipe at a circle of radius $r_0$, denoted $\Gamma_0$, around the origin of $\Omega_0$. There holds in local polar coordinates
\begin{equation}
u_0(r_0,\varphi) = \sum_{n=1}^{\infty} \left ( \mathcal{F}_n \cos(n \varphi) + \mathcal{E}_n \sin(n \varphi) \right ),
\end{equation}
cf. \cite[p.243]{Russenschuck2011}. The coefficients $\mathcal{F}_n$ and $\mathcal{E}_n$ are referred to as normal and skew coefficients, respectively. In the following, for simplicity, we assume that the symmetry of the configuration is such that all skew coefficients vanish. We consider $\mathcal{F}_n$ to be a linear functional of the solution 
\begin{equation}
\mathcal{F}_n(u_0) = \int_{\Gamma_0} \psi_n u_0 \ \mathrm{d} s = \langle \psi_n \delta_{\Gamma_0}, u_0 \rangle_{\Gamma_0}, 
\label{eq:qoi_Fourier}
\end{equation}
where $\delta_{\Gamma_0}$ and $\langle\cdot, \cdot\rangle_{\Gamma_0}$ refer to a single layer distribution associated to $\Gamma_0$ and the duality product in $\Gamma_0$, respectively. Note that in polar coordinates we have
\begin{equation}
\mathcal{F}_n(u_0) = \frac{1}{\pi} \int \limits_{0}^{2 \pi} u_0(r_0,\varphi) \cos(n \varphi) \ \mathrm{d} x.
\label{eq:multipoles}
\end{equation}
From \eqref{eq:qoi_Fourier} we see that $g_0=\psi_n \delta_{\Gamma_0}$ is the right-hand-side of the adjoint equation. Due to the single layer distribution we have 
\begin{equation}
\|\adj_0 -\adj_{h,0} \|_{L^2(\Omega_0)} =\mathcal{O}(h^{3/2}), 
\label{eq:err_adj}
\end{equation}
solely, see \cite[p.20]{becker2001optimal}.
\begin{rmk}
Collecting the Fourier coefficients as $\mathbf{f} = (\mathcal{F}_1(u_0),\mathcal{F}_2(u_0),\dots)^\top$, in magnet design the aim is often to have a harmonic distortion of 
\begin{equation}
\|\mathbf{f} - \mathbf{f}_i\|_{l^2}  \leq 10^{-4} f_i,
\end{equation}
where $\mathbf{f}_i = f_i\mathbf{e}_i$ and $\mathbf{e}_i$ denotes the $i$-th unit vector. In the case of a dipole, quadrupole and sextupole magnet, we have $i=1,2,3$, respectively, see \cite[p.242]{Russenschuck2011}.
\end{rmk}
Discrete approximations of the Fourier coefficients $\mathcal{F}_{h,n}$ are simply obtained by replacing $u_0$ with $u_{h,0}$ in \eqref{eq:multipoles} and using an error controlled adaptive numerical quadrature. To ensure a high accuracy we assume that the mesh resolves the interface as defined in \cite{hiptmair2012}. More precisely, we assume that all nodes of a triangle lie either on one side of $\Gamma_0$ or another. There holds for the FE error in the Fourier coefficients
\begin{equation}
|\mathcal{F}_n(u_0 - u_{0,h})| \leq h^{-1/2} \|u_0 - u_{0,h}\|_{L^2(\Omega_0)} = \mathcal{O}(h^{3/2}).
\label{eq:err_mp}
\end{equation}
Note that \eqref{eq:err_mp} is suboptimal, i.e., smaller than $\mathcal{O}(h^{2})$, as $\mathcal{F}_n$ is represented by a single layer distribution.

\subsection{Derivative of Magnetic Flux Density}
Magnetic deflection in a Stern-Gerlach magnet is characterized by the average partial derivative of the magnetic flux density
\begin{equation}
\mathcal{F}_{\tau}(u_0) = \frac{1}{|\Omega_0|} \int_{\Omega_0} \partial_{x_1}|\nabla u_0| \ \mathrm{d} x,
\label{eq:field_grad}
\end{equation}
\cite{DeGersem2014}, where $|\Omega_0|$ refers to the size of the domain $\Omega_0$. Contrary to the Fourier coefficients, \eqref{eq:field_grad} is a nonlinear quantity of interest with respect to the solution $u_0$. Assuming that $|\nabla u_0| >0$ in $\Omega_0$, \eqref{eq:field_grad} is well-defined as $u_0$ is smooth. During the design phase, the aim is to maximize $\mathcal{F}_{\tau}$, while minimizing the field inhomogeneity \cite{pelsoptimization}.

We emphasize that even for a higher order finite element approach the quantity $\mathcal{F}_{\tau}$ is not well-defined as $\nabla u_{h,0}$ exhibits jump discontinuities at the element interfaces. However, no difficulties arise when the cubic or quintic RBF reconstruction of the (defect corrected) FE solution is used. An estimate similar to \eqref{eq:err_mp} is beyond the scope of the paper. 

\section{Numerical Examples}
\label{sec:numerics}
Two numerical examples are given in this section to illustrate the findings. We consider an academic example on a squared domain to precisely investigate the efficiency. Then results for a Stern-Gerlach magnet are given. All results are obtained using the open-source software FEniCS \cite{Logg2012}, whereas meshes are created using Gmsh \cite{geuzaine2009gmsh}. The RBF interpolation problem \eqref{eq:rbf_system} is solved using a direct solver here.

\subsection{Academic Example}
On the domain $\Omega = [-1,1]^2$ with constant linear reluctivity and vanishing current density, an ideal sextupole ($\mathcal{F}_3$) and octupole ($\mathcal{F}_4$) field are considered. The respective solutions $u_{3} = x^3-3 x y^2$ and $u_{4} = x^4-6 x^2 y^2 + y^4$ are generated by imposing a non-homogeneous Dirichlet boundary condition on $\partial \Omega$. We use lowest order finite elements, solve the FE system using the sparse direct solver MUMPS and approximate the additional right-hand-side in \eqref{eq:dec_prim} by using a numerical Gauss quadrature of degree two. Defect correction results are provided for all polyharmonic splines presented in this paper.

\begin{table}[h!]
\small
\centering
\begin{tabular}{|c||c|c||c|c|}
\hline
$h$ & $\|u_{3} -  (\tilde{u}_{3})_h\|_{L^2(\Omega)}$ & order & $\|u_{4} -  (\tilde{u}_{4})_h\|_{L^2(\Omega)}$ & order \\
\hline
\hline
0.1414 & $1.11 \times 10^{-1}$ & - & $ 2.26\times 10^{-1}$ & -\\
\hline
0.0707 & $2.50 \times 10^{-2}$ & 2.15 & $ 5.45\times 10^{-2}$ & 2.05\\
\hline
0.0354 & $5.00 \times 10^{-3}$ & 2.32 & $ 1.13\times 10^{-2}$ & 2.27\\
\hline
0.0177 & $9.49 \times 10^{-4}$ & 2.40 & $ 2.19\times 10^{-3}$ & 2.37\\
\hline
0.0088 & $1.76 \times 10^{-4}$ & 2.43 & $ 4.15\times 10^{-4}$ & 2.40\\
\hline
\end{tabular}
\caption{Error in $L^2$-norm after primal defect correction using thin-plate splines ($k=1$).}
\label{tab:k1}
\end{table}
\begin{table}[h!]
\small
\centering
\begin{tabular}{|c||c|c||c|c|}
\hline
$h$ & $\|u_{3} -  (\tilde{u}_{3})_h\|_{L^2(\Omega)}$ & order & $\|u_{4} -  (\tilde{u}_{4})_h\|_{L^2(\Omega)}$ & order \\
\hline
\hline
0.1414 & $ 7.79\times 10^{-2}$ & - & $2.22 \times 10^{-1}$ & -\\
\hline
0.0707 & $ 1.25\times 10^{-2}$ & 2.64 & $3.62 \times 10^{-2}$ & 2.61\\
\hline
0.0354 & $ 1.78\times 10^{-3}$ & 2.81 & $5.28 \times 10^{-3}$ & 2.78\\
\hline
0.0177 & $ 2.39\times 10^{-4}$ & 2.89 & $7.16 \times 10^{-4}$ & 2.88\\
\hline
0.0088 & $ 3.13\times 10^{-5}$ & 2.94 & $9.36 \times 10^{-5}$ & 2.94\\
\hline
\end{tabular}
\caption{Error in $L^2$-norm after primal defect correction using cubics ($k=2$).}
\label{tab:k2}
\end{table}
\begin{table}[h!]
\small
\centering
\begin{tabular}{|c||c|c||c|c|}
\hline
$h$ & $\|u_{3} -  (\tilde{u}_{3})_h\|_{L^2(\Omega)}$ & order & $\|u_{4} -  (\tilde{u}_{4})_h\|_{L^2(\Omega)}$ & order \\
\hline
\hline
0.1414 & $1.47 \times 10^{-2}$ & - & $1.48 \times 10^{-1}$ & -\\
\hline
0.0707 & $1.14 \times 10^{-3}$ & 3.68 & $1.11 \times 10^{-2}$ & 3.74\\
\hline
0.0354 & $7.78 \times 10^{-5}$ & 3.88 & $8.07 \times 10^{-4}$ & 3.78\\
\hline
0.0177 & $5.15 \times 10^{-6}$ & 3.92 & $5.44 \times 10^{-5}$ & 3.89\\
\hline
0.0088 & $3.36 \times 10^{-7}$ & 3.94 & $3.54 \times 10^{-6}$ & 3.94\\
\hline
\end{tabular}
\caption{Error in $L^2$-norm after primal defect correction using quintics ($k=3$).}
\label{tab:k3}
\end{table}
To accurately determine the convergence order of the defect correction scheme we consider a structured grid and the compute the error in the $L^2$-norm. Tables \ref{tab:k1},\ref{tab:k2},\ref{tab:k3} show the results for thin-plate splines, cubics and quintics, respectively. The predicted convergence rates \eqref{eq:rbf_order} are well-observed.

\begin{figure}[t!]
\begin{minipage}[t!]{0.5\textwidth}
\includegraphics[width=\textwidth,trim=0cm 0cm 5.6cm 0cm,clip]{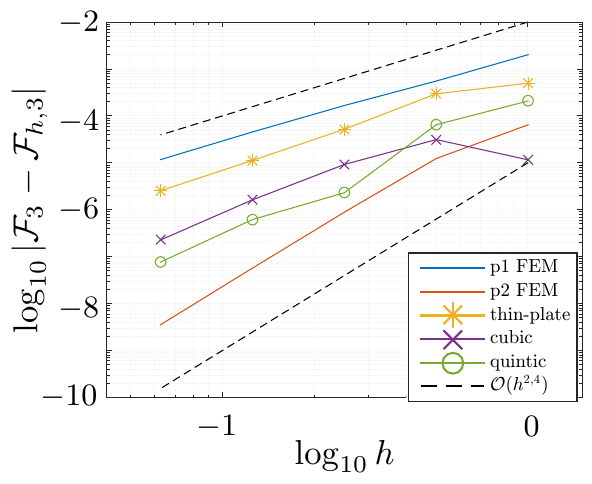}
\end{minipage}
\begin{minipage}[t!]{0.5\textwidth}
\includegraphics[width=\textwidth,trim=0cm 0cm 5.6cm 0cm,clip]{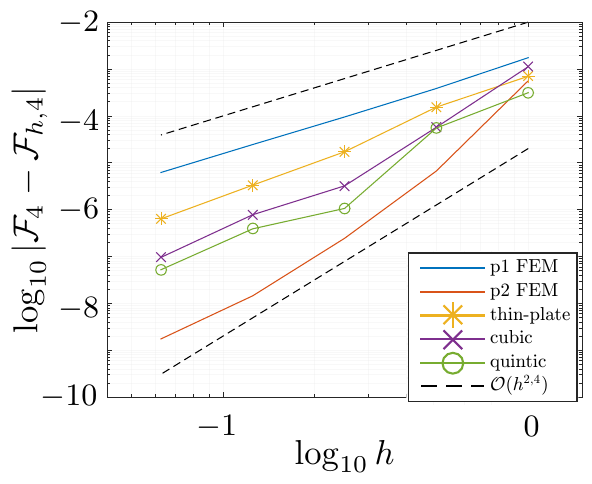}
\end{minipage}
\caption{Discretization error in Fourier coefficients for standard finite element approaches and primal defect correction using polyharmonic splines. Left: sextupole component $\mathcal{F}_3$. Right: octopole component $\mathcal{F}_4$.}
\label{fig:dc_prim}
\end{figure}

In order to extract Fourier coefficients we employ an unstructured mesh of maximum mesh size $h$, which is aligned at the reference circle $\Gamma_0$. Several steps of uniform mesh refinement are carried out. Concerning the Fourier coefficients, at a reference radius $r_0 = 0.2$ we compute $\mathcal{F}_3=0.008$ and $\mathcal{F}_4=0.0016$ in the case of $u_{3}$ and $u_{4}$, respectively. These coefficients are computed using an adaptive quadrature with an absolute error smaller than $10^{-15}$. In Figure \ref{fig:dc_prim}, the errors in the Fourier coefficients of standard linear and quadratic finite elements are depicted and compared to primal defect correction using polyharmonic splines. Again, a significant improvement of the convergence rate due to defect correction can be observed, although the convergence orders are more difficult to extract compared to the $L^2$ norm.

\begin{figure}[t!]
\begin{minipage}[t!]{0.5\textwidth}
\includegraphics[width=\textwidth,trim=0cm 0cm 5.6cm 0cm,clip]{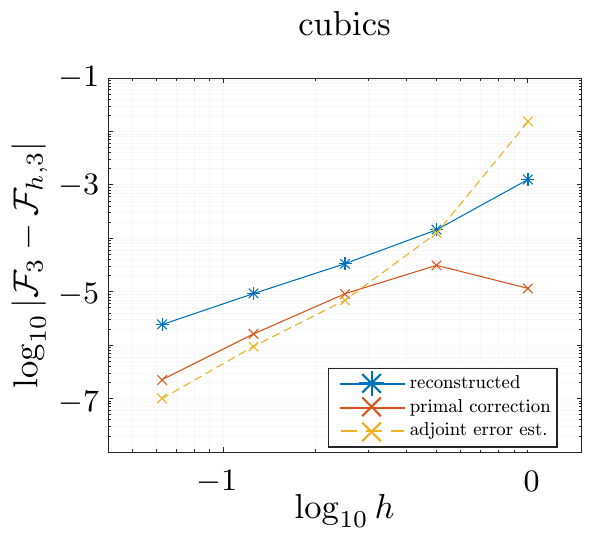}
\end{minipage}
\begin{minipage}[t!]{0.5\textwidth}
\includegraphics[width=\textwidth,trim=0cm 0cm 5.6cm 0cm,clip]{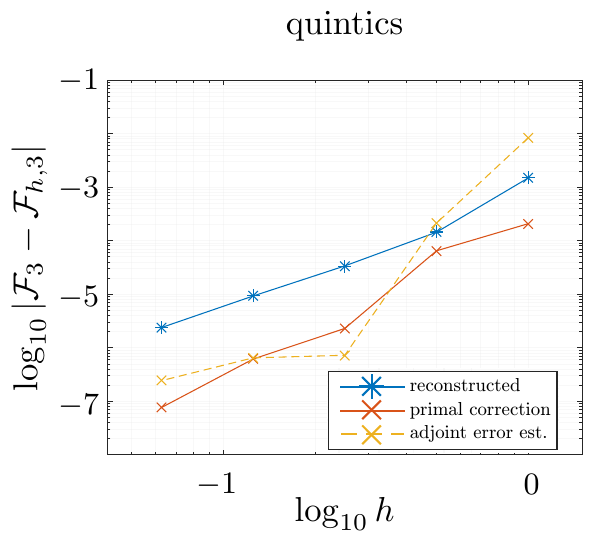}
\end{minipage}
\caption{Discretization error in Fourier coefficient $\mathcal{F}_3$ for reconstructed defect corrected solution. The remaining error after defect correction is estimated using adjoint correction. Left: cubics ($k=2$). Right: quintics ($k=3$).}
\label{fig:dc_prim_2}
\end{figure}

In practice no reference solution is available and the remaining error needs to be estimated. To this end we employ adjoint correction as outlined in Section \ref{sec:dec}. As it turns out that it is difficult to assemble the right-hand-side of the adjoint equation, due to the single layer distribution. This problem is circumvented here, by finding a volume based formulation of the quantity of interest in the interior of $\Gamma_0$ using the divergence theorem. In Figure \ref{fig:dc_prim_2} the errors in the Fourier coefficients are depicted for the reconstructed solution and the solution after primal defect correction. It can be observed, that reconstructing the solution solely, does not improve the convergence order. Hence, the interest in defect correction. The remaining error is estimated using the adjoint approach. We observe, that the error estimator is accurate for finer meshes but overestimates the true error for coarse meshes. 
~\newline

\subsection{Stern-Gerlach Magnet}
We consider the example of a Rabi-type Stern-Gerlach magnet. Details on geometry and the numerical setup can be found in \cite{masschaele2012design,pelsoptimization}. It should be noted that we consider a linear material here with reluctivity $\nu = 1/(\mu_0 \mu_{\mathrm{r}})$, where $\mu_0=4 \pi 10^{-7}$ H/m and $\mu_{\mathrm{r}} = 1000$, respectively. This explains deviations in the results in the order of $2 \%$ compared to \cite{pelsoptimization}, where a nonlinear material was considered. A constant current of $2600$ A is imposed. In this example the FE system is solved using the conjugate gradient algorithm with algebraic multigrid preconditioning. Pollution error control is achieved using adaptive mesh refinement based on the error indicator $\eta$ given in Section \ref{sec:dec}. A global-adaptive algorithm, as given in \cite[Section 3.2]{babuvska1997pollution}, is used. It consists in refining all elements with local indicator $\eta_K \geq \gamma \max_{T \in \mathcal{T}_h} \eta_T$. For this example $\gamma=0.5$ is chosen and the process is stopped after $n_{\mathrm{ref}}$ iterations. 

\begin{figure}[t!]
\centering
\includegraphics[width = 0.7\textwidth]{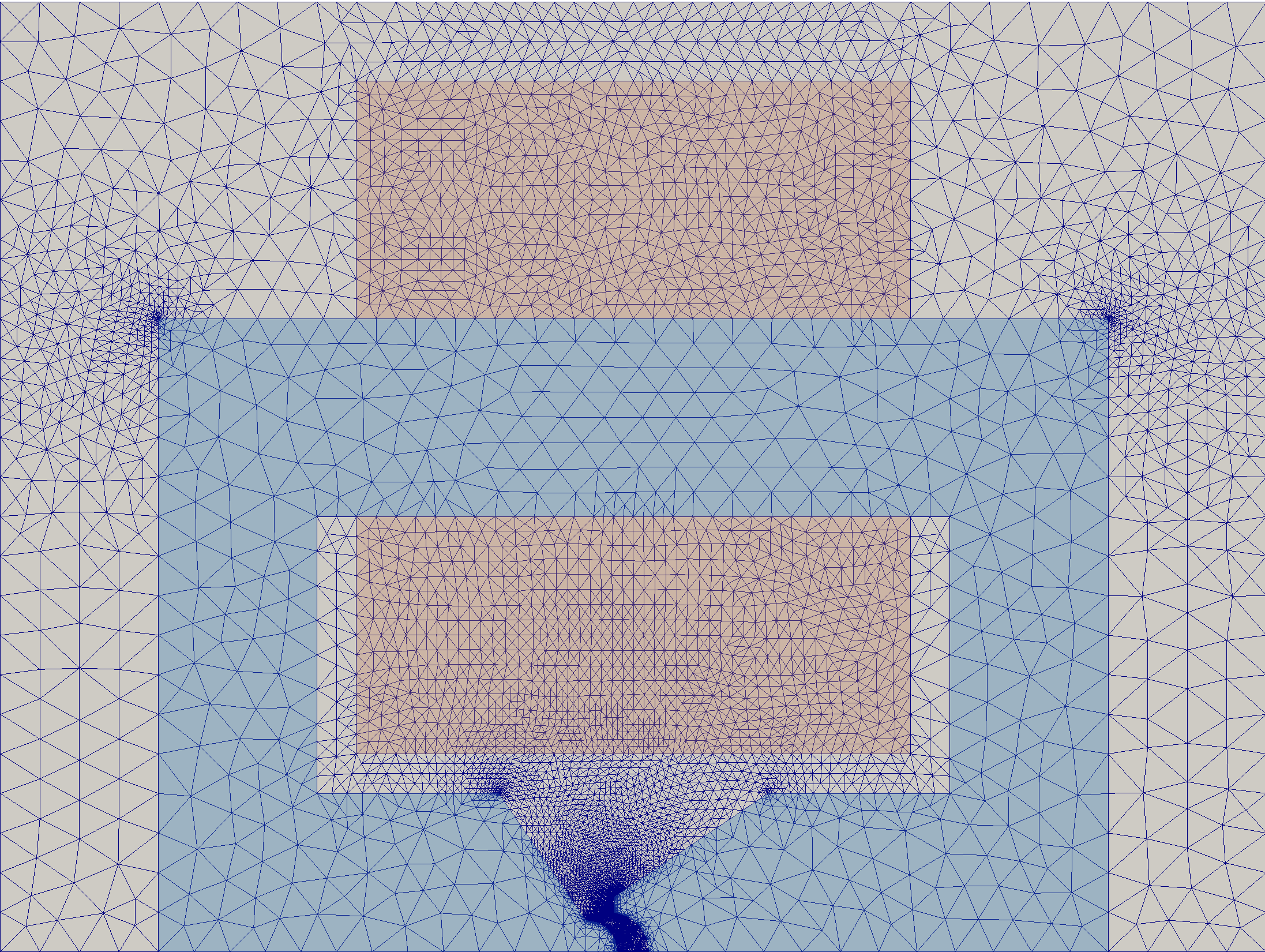}
\caption{Geometry and adaptively refined mesh of Rabi-type Stern-Gerlach magnet. Domains $\Omega_\con$ and $\Omega_\iso$ in blue and red/grey color, respectively.}
\label{fig:sg}
\end{figure}

The geometry and an adaptively refined mesh with $\eta_{\mathrm{rel}}=0.09$, where $\eta_{\mathrm{rel}}$ refers to $\eta$ divided by the $H^1$-norm of the solution, are depicted in Figure \ref{fig:sg}. It can be clearly observed how the mesh is refined at the re-entrant corners. Furthermore, although refinement using the global quantity $\eta_{\mathrm{rel}}$ does not guarantee sufficient refinement in $\Omega_0$, the mesh in the beam area is very dense due to the singularities and the complicated geometry of the surrounding iron yoke.

\begin{figure}[t!]
\centering
\begin{minipage}{0.3\textwidth}
\includegraphics[width=\textwidth]{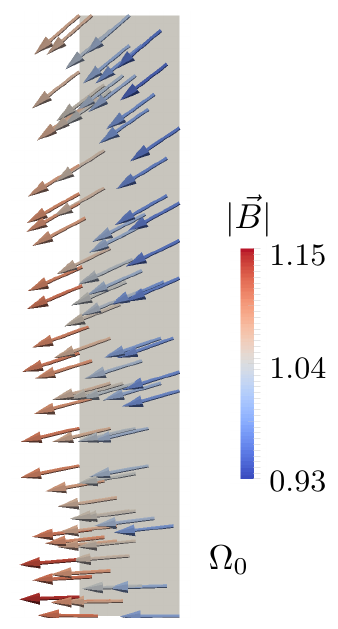}
\end{minipage}
\begin{minipage}{0.3\textwidth}
\includegraphics[width=\textwidth]{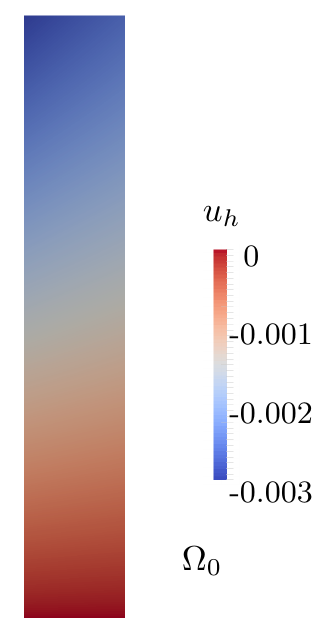}
\end{minipage}
\begin{minipage}{0.3\textwidth}
\includegraphics[width=\textwidth]{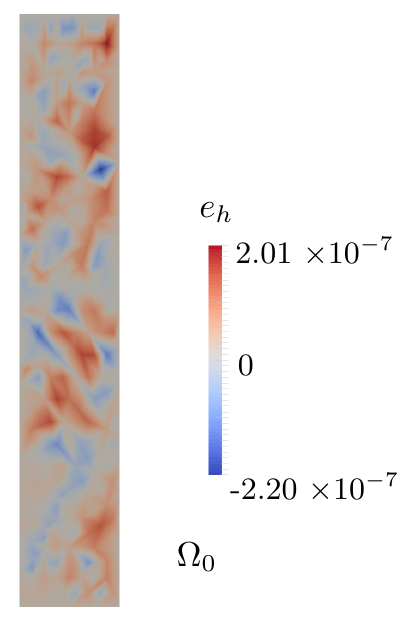}
\end{minipage}
\caption{From left to right: magnetic flux density, potential and estimated error before defect correction in the region $\Omega_0$, inside the air gap, where the quantity of interest is evaluated.}
\label{fig:sg_quantities}
\end{figure}

Figure \ref{fig:sg_quantities} depicts the potential and the magnetic flux density in $\Omega_0$, the region inside of the airgap, where the quantity of interest is evaluated. Additionally, the solution of \eqref{eq:dec_prim_loc} is plotted, which is an estimate of the local error before defect correction. Numerical results for the quantity of interest are given in Table \ref{tab:sg_fine} for different adaptively refined meshes. Defect correction is carried out twice using quintics ($k=3$), yielding the most accurate results, as shown in the previous example. As a reference solution a second degree FE solution is employed. The quantity of interest is evaluated by projecting $\nabla u_h$ onto the space of continuous vector functions of degree one. Hence, the (weak) derivative of $\nabla u_h$ exists and $\mathcal{F}_{\tau}$ is well-defined. As seen from Table \ref{tab:sg_fine}, defect correction combined with globally-adaptive mesh refinement yields an average uni-directional derivative of the magnetic flux density of $-247.501$ T/m with an error of below $1 \%$ with respect to the reference solution. Also, a global error indicator of $\eta_{\mathrm{rel}} \approx 0.1$ seems to be sufficient to ensure a pollution error with the same order of magnitude. It should be noted that $\eta_{\mathrm{rel}}$ typically overestimates the true error, as unknown constants are neglected. 
\begin{table}[h!]
\small
\centering
\begin{tabular}{|c||c|c||c|c|}
\hline
 & \multicolumn{2}{|c||}{defect correction ($k=3$)} & \multicolumn{2}{|c|}{higher order reference} \\
\hline
$n_{\mathrm{ref}}$ & $\eta_{\mathrm{rel}}$ & $\mathcal{F}_{\tau}$ T/m & $\eta_{\mathrm{rel}}$ &  $\mathcal{F}_{\tau}$ T/m\\
\hline
10 & 0.1974 & -244.878 & 0.0856 & -246.866 \\
\hline
20 & 0.0581 & -246.273 & 0.0165 & -246.689 \\
\hline
30 & 0.0273 & -246.806 & $<0.0001$ & -246.813 \\
\hline
\end{tabular}
\caption{Error in $\mathcal{F}_\tau$ using two times defect correction with quintics ($k=3$) and a higher order reference solution with post-processing for a different number of adaptive refinement levels.}
\label{tab:sg_fine}
\end{table}


\subsection{Discussion of Costs}
The quadratic FE approach seems to outperform the first order FE in combination with polyharmonic spline defect correction. However, we do not view these methods as competitors, as defect correction can be applied in combination with higher order FE as well. Moreover, as stressed in the introduction, the aim is also to improve the differentiability of the solution. Yet, we briefly compare the complexity, to give a better impression of the associated costs: a state-of-the art quadratic FE method has $\mathcal{O}(N)$ complexity and only one linear system needs to be solved. This is also true for a linear FE method. The additional efforts for defect correction are 
\begin{itemize}
\item
the solution of \eqref{eq:rbf_system} with $\mathcal{O}(N^3)$ operations ($\mathcal{O}(N \log N)$ with the fast multipole method),
\item
the assembly of the right-hand-side in \eqref{eq:dec_prim} with $\mathcal{O}(N^2)$ operations ($\mathcal{O}(N)$ with the partition of unity method \cite{Wendland02fastevaluation}),
\item the solution of the error equation \eqref{eq:dec_prim} with complexity $\mathcal{O}(N)$,
\item RBF reconstruction of the error by solving \eqref{eq:rbf_system} with $\mathcal{O}(N^3)$ ($\mathcal{O}(N \log N)$) operations.
\end{itemize}
For the local defect correction method the costs with respect to the local and global number of degrees of freedom $N_0$ and $N$ are: solution of \eqref{eq:rbf_system} ($\mathcal{O}(N_0^3)$ or $\mathcal{O}(N_0 \log N_0)$), assembly of the right-hand-side in \eqref{eq:dec_prim} ($\mathcal{O}(N N_0)$ or $\mathcal{O}(N)$), solution of \eqref{eq:dec_prim} ($\mathcal{O}(N)$), solution of \eqref{eq:rbf_system} ($\mathcal{O}(N_0^3)$ or $\mathcal{O}(N_0 \log N_0)$) operations.

\section{Conclusion}
In this work, a defect correction scheme for the accurate numerical approximation of magnetic fields was presented. The post-processing was achieved using radial basis functions and is general as it allows for unstructured grids. Using adjoint techniques, the remaining error after defect correction could be estimated. It was outlined how defect correction can be used to improve the local approximation error in the case of interface problems, whereas the pollution error was addressed using explicit residual error estimators. Convergence estimates for the defect correction scheme were discussed and observed in numerical examples. Finally, accurate simulations results were obtained for a two-dimensional model of an actually existing Stern-Gerlach magnet.

\section*{Acknowledgment}
This work was supported  by the `Excellence Initiative' of German Federal and State Governments and the Graduate School CE at Technische Universit\"at Darmstadt, by the Deutsche Forschungsgemeinschaft under SFB 634 and by the project nanoCOPS founded by the European Union.


\section*{References}

\begin{thebibliography}{10}
\expandafter\ifx\csname url\endcsname\relax
  \def\url#1{\texttt{#1}}\fi
\expandafter\ifx\csname urlprefix\endcsname\relax\def\urlprefix{URL }\fi
\expandafter\ifx\csname href\endcsname\relax
  \def\href#1#2{#2} \def\path#1{#1}\fi

\bibitem{Russenschuck2011}
S.~Russenschuck, Field computation for accelerator magnets: analytical and
  numerical methods for electromagnetic design and optimization, John Wiley and
  Sons, 2011.

\bibitem{salon1995finite}
S.~J. Salon, Finite element analysis of electrical machines, Vol. 101, Kluwer
  academic publishers, Boston, USA, 1995.

\bibitem{Kurz1998}
S.~Kurz, W.~M. Rucker, J.~Fetzer, Coupled {BEM}-{FEM} methods for 3{D} field
  calculations with iron saturation, Tech. rep., Proceedings of the First
  International ROXIE users meeting and workshop, CERN (1998).

\bibitem{DeGersem2002}
H.~{De~Gersem}, M.~Clemens, T.~Weiland, Coupled finite-element,
  spectral-element discretisation for models with circular inclusions and
  far-field domains, IEE Proceedings-Science, Measurement and Technology 149
  (2002) 237--241.

\bibitem{demkowicz2006computing}
L.~Demkowicz, Computing with hp-Adaptive Finite Elements: Volume 1 One and Two
  Dimensional Elliptic and {M}axwell problems, Chapman and Hall/CRC, Boca
  Raton, USA, 2006.

\bibitem{dular2009posteriori}
P.~Dular, A posteriori error estimation of finite element solutions via the
  direct use of higher order hierarchal test functions, IEEE transactions on
  magnetics 45~(3) (2009) 1360--1363.

\bibitem{DeGersem2014}
H.~{De~Gersem}, B.~Masschaele, T.~Roggen, E.~Janssens, N.~Tung, Improved field
  post-processing for a {S}tern--{G}erlach magnetic deflection magnet,
  International Journal of Numerical Modelling: Electronic Networks, Devices
  and Fields 27 (2014) 472--484.

\bibitem{pelsoptimization}
A.~Pels, Z.~Bontinck, J.~Corno, H.~{De~Gersem}, S.~Sch\"{o}ps, Optimization of
  a stern-gerlach magnet by magnetic field-circuit coupling and isogeometric
  analysis.

\bibitem{Stetter1978}
H.~J. Stetter, The defect correction principle and discretization methods,
  Numerische Mathematik 29 (1978) 425--443.

\bibitem{bohmer2012defect}
K.~B{\"o}hmer, H.~J. Stetter, Defect correction methods: theory and
  applications, Vol.~5, Springer Science \& Business Media, 2012.

\bibitem{Giles2002}
M.~B. Giles, E.~S{\"u}li, Adjoint methods for {PDE}s: a posteriori error
  analysis and postprocessing by duality, Acta Numerica 11 (2002) 145--236.

\bibitem{pierce2000adjoint}
N.~A. Pierce, M.~Giles, Adjoint recovery of superconvergent functionals from
  pde approximations, SIAM {R}eview 42~(2) (2000) 247--264.

\bibitem{Pierce2004}
N.~A. Pierce, M.~B. Giles, Adjoint and defect error bounding and correction for
  functional estimates, Journal of Computational Physics 200 (2004) 769--794.

\bibitem{basumatary2014defect}
M.~Basumatary, G.~Natarajan, S.~C. Mishra, Defect correction based velocity
  reconstruction for physically consistent simulations of non-{N}ewtonian flows
  on unstructured grids, Journal of Computational Physics 272 (2014) 227--244.

\bibitem{bonaventura2011kernel}
L.~Bonaventura, A.~Iske, E.~Miglio, Kernel-based vector field reconstruction in
  computational fluid dynamic models, International Journal for Numerical
  Methods in Fluids 66~(6) (2011) 714--729.

\bibitem{cullity2011introduction}
B.~D. Cullity, C.~D. Graham, Introduction to magnetic materials, John Wiley \&
  Sons, 2011.

\bibitem{Pierce2000}
N.~A. Pierce, M.~B. Giles, Adjoint recovery of superconvergent functionals from
  pde approximations, SIAM Review 42 (2000) 247--264.

\bibitem{barrett1988optimal}
J.~Barrett, G.~Moore, K.~Morton, Optimal recovery in the finite-element method,
  part 2: Defect correction for ordinary differential equations, IMA {J}ournal
  of {N}umerical {A}nalysis 8~(4) (1988) 527--540.

\bibitem{buhmann_radial_2000}
M.~D. Buhmann, Radial basis functions, Acta Numerica 2000 9 (2000) 1--38.

\bibitem{gumerov2007fast}
N.~A. Gumerov, R.~Duraiswami, Fast radial basis function interpolation via
  preconditioned {K}rylov iteration, SIAM Journal on Scientific Computing
  29~(5) (2007) 1876--1899.

\bibitem{johnson2001}
M.~Johnson, The {L}2-approximation order of surface spline interpolation,
  Mathematics of {C}omputation 70~(234) (2001) 719--737.

\bibitem{wendland_meshless_1999}
H.~Wendland, Meshless {G}alerkin methods using radial basis functions,
  Mathematics of Computation of the American Mathematical Society 68~(228)
  (1999) 1521--1531.

\bibitem{fuselier2007refined}
E.~J. Fuselier~Jr, Refined error estimates for matrix-valued radial basis
  functions, Ph.D. thesis, Texas A\&M University (2007).

\bibitem{giles2004progress}
M.~B. Giles, N.~Pierce, E.~S{\"u}li, Progress in adjoint error correction for
  integral functionals, Computing and Visualization in Science 6~(2-3) (2004)
  113--121.

\bibitem{hackbusch1984local}
W.~Hackbusch, Local defect correction method and domain decomposition
  techniques, in: Defect correction methods, Springer, 1984, pp. 89--113.

\bibitem{babuvska1997pollution}
I.~Babu{\v{s}}ka, T.~Strouboulis, S.~Gangaraj, C.~Upadhyay, Pollution error in
  the h-version of the finite element method and the local quality of the
  recovered derivatives, Computer Methods in Applied Mechanics and Engineering
  140~(1) (1997) 1--37.

\bibitem{oden1996local}
J.~T. Oden, Y.~Feng, Local and pollution error estimation for finite element
  approximations of elliptic boundary value problems, Journal of computational
  and Applied Mathematics 74~(1) (1996) 245--293.

\bibitem{ainsworth1997posteriori}
M.~Ainsworth, J.~T. Oden, A posteriori error estimation in finite element
  analysis, Computer Methods in Applied Mechanics and Engineering 142~(1)
  (1997) 1--88.

\bibitem{becker2001optimal}
R.~Becker, R.~Rannacher, An optimal control approach to a posteriori error
  estimation in finite element methods, Acta Numerica 2001 10 (2001) 1--102.

\bibitem{hiptmair2012}
R.~Hiptmair, J.~Li, J.~Zou, Convergence analysis of finite element methods for
  $\mathrm{{H}}(\mathrm{curl};{\Omega})$-elliptic interface problems,
  Numerische Mathematik 122~(3) (2012) 557--578.

\bibitem{Logg2012}
A.~Logg, K.~A. Mardal, G.~N. Wells, Automated Solution of Differential
  Equations by the Finite Element Method, Springer, 2012.

\bibitem{geuzaine2009gmsh}
C.~Geuzaine, J.~Remacle, Gmsh: A 3-{D} finite element mesh generator with
  built-in pre-and post-processing facilities, International Journal for
  Numerical Methods in Engineering 79~(11) (2009) 1309--1331.

\bibitem{masschaele2012design}
B.~Masschaele, T.~Roggen, H.~{De~Gersem}, W.~Janssens, T.~T. Nguyen, Design of
  a strong gradient magnet for the deflection of nanoclusters, Applied
  Superconductivity, IEEE Transactions on 22~(3) (2012) 3700604--3700604.

\bibitem{Wendland02fastevaluation}
H.~Wendland, Fast evaluation of radial basis functions: Methods based on
  partition of unity, in: Approximation Theory X: Wavelets, Splines, and
  Applications, Vanderbilt University Press, 2002, pp. 473--483.

\end{thebibliography}
\end{document}